\documentclass[11pt,reqno]{amsart}
\usepackage[T1]{fontenc}
\usepackage[letterpaper,textwidth=5.45in,textheight=8.4in,centering]{geometry}
\usepackage{lmodern,amssymb,mathtools}
\usepackage[hidelinks]{hyperref}
\newtheorem{theorem}{Theorem}[section]
\newtheorem{lemma}[theorem]{Lemma}
\newtheorem{corollary}[theorem]{Corollary}
\theoremstyle{remark}
\newtheorem{example}[theorem]{Example}
\newcommand{\Inn}{\operatorname{Inn}}
\newcommand{\Aut}{\operatorname{Aut}}
\newcommand{\supp}{\operatorname{supp}}
\newcommand{\Rad}{\operatorname{Rad}}
\newcommand{\RP}{\mathbb{R}P}
\newcommand{\Z}{\mathbb Z}
\newcommand{\F}{\mathbb F}
\newcommand{\PP}{\mathcal P}
\newcommand{\SR}{\mathcal R}

\title[Homology of subrack lattices]{Torsion in the homology of subrack lattices}
\author{Amir Saki}
\address{Universit\'e du Qu\'ebec \`a Trois-Rivi\`eres, Trois-Rivi\`eres, Qu\'ebec, Canada}
\email{amir.saki.math@gmail.com}
\author{Usef Faghihi}
\address{Universit\'e du Qu\'ebec \`a Trois-Rivi\`eres, Trois-Rivi\`eres, Qu\'ebec, Canada}
\email{usef.faghihi@uqtr.ca}

\begin{document}
\begin{abstract}
Heckenberger, Shareshian, and Welker asked whether the order complex of the nonempty proper subracks of a rack can have torsion in its integral homology. We construct a finite quandle for which the homology in degree \(298\) has a direct summand isomorphic to \(\Z/2\Z\). More generally, for each finite simplicial complex \(K\) with at least one vertex, we construct a finite conjugation quandle whose proper subrack order complex contains a suspension of \(K\) as a homotopy retract. It follows that every finite abelian group occurs as a direct summand in the integral homology of such an order complex. Heckenberger, Shareshian, and Welker also asked whether a proper subrack order complex can be acyclic. For finite racks, we show that this is equivalent to asking whether the order complex is contractible for some union of conjugacy classes in a finite group with trivial solvable radical.
\end{abstract}
\maketitle

\section{Introduction}

Let \(R\) be a finite rack. We write \(\SR(R)\) for the lattice of subracks of \(R\), ordered by inclusion, and \(\Delta_R\) for the simplicial complex whose simplices are chains of nonempty proper subracks. We call a complex acyclic if its realization is nonempty and all its reduced integral homology groups vanish.

In \cite[Section~5]{HSW}, Heckenberger, Shareshian, and Welker asked whether the integral homology of \(\Delta_R\) can have torsion. They also asked whether \(\Delta_R\) can be acyclic. We consider the acyclicity question for finite racks. For the conjugation rack of an entire finite group, the authors showed that \(\Delta_R\) is homotopy equivalent to a sphere \cite{HSW}. Saki and Kiani showed that the subrack lattice of every finite rack is complemented \cite{SK}. Kayacan related the homotopy type of \(\Delta_R\) for decomposable conjugation racks to the order complex of the poset of parabolic subracks \cite{Kayacan}.

We answer the torsion question affirmatively. Given a finite simplicial complex \(K\), we choose conjugacy classes in a finite group so that a suspension of \(K\) is a homotopy retract of the order complex of the resulting quandle. Applied to the six-vertex triangulation of the real projective plane, the construction gives
\[
\Z/2\Z\ \text{as a direct summand of }\
\widetilde H_{298}(\Delta_R;\Z).
\]
Applied to finite simplicial models of Moore spaces, it gives every finite abelian group as a direct summand.

The existence of a finite rack with acyclic \(\Delta_R\) remains open. We show that it suffices to consider unions of conjugacy classes that generate a finite group with trivial solvable radical. Within this class, the existence of an acyclic example is equivalent to the existence of a contractible example.

\section{The subgroup poset associated with a rack}

A rack is a set \(R\) with a binary operation \(x\triangleright y\) such that each left translation \(L_x\) is a bijection and
\(x\triangleright(y\triangleright z)=(x\triangleright y)\triangleright(x\triangleright z)\)
for all \(x,y,z\in R\). A rack is a quandle if \(x\triangleright x=x\) for each \(x\in R\). Let \(\Inn(R)\) denote the group generated by the left translations. A subrack is a subset on which the restricted operation defines a rack. If \(G\) is a group and \(R\subseteq G\) is a union of conjugacy classes, then \(R\) is a quandle under \(x\triangleright y=xyx^{-1}\). We call \(R\) a conjugation quandle.

Let \(G\) be a finite group and let \(C_1,\dots,C_r\) be a list of conjugacy classes whose union generates \(G\). We allow repetitions in the list. Put
\[
\PP(G;C_1,\dots,C_r)
=\{H<G:H\cap C_i\ne\varnothing\text{ for each }i\}.
\]
Thus \(\PP(G;C_1,\dots,C_r)\) is the poset of proper subgroups of \(G\) that meet every \(C_i\). Repeating a class gives the same poset. For each poset \(P\), let \(\Delta(P)\) denote the simplicial complex of chains in \(P\). We write \(*\) for the simplicial join. The order complex of the empty poset consists of the empty face. We denote it by \(S^{-1}\), so \(S^{-1}*K=K\).

For conjugation racks, Kayacan's decomposition expresses the order complex as a join of a sphere and the order complex of the poset of parabolic subracks \cite{Kayacan}. The next lemma gives a subgroup formulation and applies to any finite rack through its inner automorphism group.

\begin{lemma}\label{lem:reduction}
Let \(R\) be a finite nonempty rack with \(r\) orbits under \(G=\Inn(R)\). For each \(i\), let \(C_i\) be the image of the \(i\)-th orbit under \(x\mapsto L_x\). Then
\[
\Delta_R \simeq S^{r-2} * \Delta\bigl(\PP(G;C_1,\dots,C_r)\bigr).
\]
The same formula holds when \(R\) is the union of \(r\) distinct conjugacy classes \(C_i\) generating a finite group \(G\), with the conjugation operation.
\end{lemma}

\begin{proof}
We give a proof that works for arbitrary finite racks. In the first case set \(\lambda(x)=L_x\). In the second case set \(\lambda(x)=x\) and let \(G\) act on \(R\) by conjugation. In both cases \(\lambda(R)\) generates \(G\), and \(\lambda(x)\) acts as \(L_x\) for each \(x\in R\). Also, \(\lambda(gx)=g\lambda(x)g^{-1}\) for all \(g\in G\) and \(x\in R\). For the action of \(\Inn(R)\), this follows from \(L_yL_xL_y^{-1}=L_{y\triangleright x}\). In particular, each orbit has a conjugacy class as its image.

Let \(F\) be the poset of proper subracks meeting every \(G\)-orbit. Write \(P\) for the subgroup poset in the statement. Define \(F\to P\) by \(A\mapsto\langle\lambda(A)\rangle\) and \(P\to F\) by \(H\mapsto\lambda^{-1}(H)\). To see that these maps are well defined, observe that a subrack \(A\) is invariant under the action of \(\langle\lambda(A)\rangle\), so a subrack in \(F\) cannot generate \(G\). Also, \(\lambda^{-1}(H)\) is a proper subrack because \(\lambda(R)\) generates \(G\). For each \(A\in F\), we have \(A\subseteq\lambda^{-1}(\langle\lambda(A)\rangle)\). For each \(H\in P\), we have \(\langle\lambda(\lambda^{-1}(H))\rangle\leq H\). These inequalities give homotopies between the composites and the identity maps by the order homotopy theorem \cite[Theorem~10.11]{Bjorner}.

Let \(B_r\) be the poset of nonempty proper subsets of \(\{1,\dots,r\}\). Order the disjoint union \(B_r\sqcup F\) by retaining the two given orders and putting every member of \(B_r\) below every member of \(F\). This is the ordinal sum \(B_r\oplus F\). Map the nonempty proper subracks of \(R\) to this poset as follows. Send a subrack missing an orbit to the set of orbit indices it meets and send each \(A\in F\) to itself. The inverse image of the lower ideal at a member of \(B_r\) has a greatest element, the union of the corresponding orbits. For each \(B\in F\), consider the inverse image of the lower ideal at \(B\). For each \(A\) in that inverse image, let \(R_A\) be the union of the orbits met by \(A\). The order-preserving maps given by
\[
 A\ \leq\ A\vee(B\cap R_A)\ \geq\ B\cap R_A\ \leq\ B
\]
take values in that inverse image. If \(A\) misses an orbit, the join is contained in the proper subrack \(R_A\). If \(A\) meets all orbits, then \(A\leq B\) and the join is \(B\). Since \(B\) meets each orbit, \(B\cap R_A\) is nonempty. Hence this inverse image is contractible. Quillen's fiber lemma \cite[Proposition~1.6]{Quillen} gives \(\Delta_R\simeq\Delta(B_r)*\Delta(F)\). Since \(\Delta(B_r)\) is the subdivision of the boundary of an \((r-1)\)-simplex, the assertion follows.
\end{proof}

\section{A construction for finite simplicial complexes}\label{sec:construction}

\begin{theorem}\label{thm:realization}
Let \(K\) be a finite simplicial complex with at least one vertex. There is a finite group \(G\) and a union \(R\) of \(r\) conjugacy classes generating \(G\) such that \(\Sigma^{r-1}|K|\) is a homotopy retract of \(|\Delta_R|\). Here \(\Sigma^j\) denotes \(j\) successive suspensions. In particular, for each \(i\geq0\), the group \(\widetilde H_i(K;\Z)\) is isomorphic to a direct summand of \(\widetilde H_{i+r-1}(\Delta_R;\Z)\).
\end{theorem}

\begin{proof}
Let \(\{1,\dots,m\}\) be the vertex set of \(K\). A minimal nonface is a set of vertices that is not a face of \(K\) and whose proper subsets are all faces. Choose an odd prime \(\ell\) at least as large as the size of each minimal nonface. If every subset is a face, choose any odd prime. Put \(V=\F_\ell^m\). The support of a vector is the set of indices of its nonzero coordinates. The nonzero vectors on a line \(L\leq V\) all have the same support. We call \(L\) forbidden if this support is a nonface of \(K\). For each subspace \(U\) of \(V\), let \(\supp(U)\) be the union of the supports of its vectors. Then
\begin{equation}\label{eq:support}
U\text{ contains no forbidden line}\quad\Longleftrightarrow\quad \supp(U)\in K.
\end{equation}
Indeed, assume that \(\supp(U)\) contains a minimal nonface \(E\). For each \(i\in E\), the kernel of the \(i\)-th coordinate on \(U\) is a hyperplane. If \(d=\dim U\), the union of these hyperplanes has at most \(1+|E|(\ell^{d-1}-1)<\ell^d\) elements. A vector outside the union spans a forbidden line. The converse follows because every subset of a face is a face.

Let \(T=A_{\ell+2}\), and fix an \(\ell\)-cycle \(s\in T\). The group \(T\) is nonabelian simple. Also, \(\Aut(T)\cong S_{\ell+2}\), since \(\ell+2\geq5\) and \(\ell+2\ne6\) \cite{DM}. Since \(\ell+2<\ell^2\), no permutation in \(S_{\ell+2}\) has order \(\ell^2\). Let \(N=T^V\) be the group of functions \(V\to T\) with pointwise multiplication. Let \(V\) act on \(N\) by translation, so that
\[
(q\cdot f)(v)=f(v-q)
\qquad(q,v\in V,\ f\in N).
\]
Put \(G=N\rtimes V\) for this action. Its quotient map to \(V\) sends \((f,q)\) to \(q\). Define \(R\) to be the union of the following conjugacy classes of \(G\).
\begin{enumerate}
\item For each nonidentity conjugacy class \(C\) of \(T\), include the base class \(C^V\).
\item For each standard basis vector \(e_j\) of \(V\), include the class of \((1,e_j)\).
\item For each forbidden line, choose a nonzero vector \(q\) spanning it. On each orbit of translation by \(q\) in \(V\), choose one point. Let \(f_q\) have value \(s\) at the chosen points and value \(1\) elsewhere. Include the class of \((f_q,q)\).
\end{enumerate}
Each set \(C^V\) in (1) is a single class, since conjugation by \(N\) can be performed independently in each coordinate. The classes in (1) generate \(N\). In fact, for each coordinate, the subgroup generated by those classes contains a nonidentity element supported in that coordinate, obtained as a quotient of two functions of a fixed class which differ only there. This subgroup is normal in \(G\), so simplicity gives the entire copy of \(T\) in that coordinate. The classes in (2) generate the quotient \(V\). Thus \(R\) generates \(G\).

The first two types of classes control how a subgroup projects onto the coordinates of \(N\) and onto \(V\). The third type encodes the nonfaces of \(K\).

Write \(P=\PP(G;C_1,\dots,C_r)\), where \(C_1,\dots,C_r\) are the selected classes. Let \(H\in P\), and set \(D=H\cap N\). By (2), \(H\) maps onto \(V\). By (1), the projection of \(D\) onto each coordinate meets every conjugacy class of \(T\), and hence equals \(T\). Indeed, the union of the conjugates of a proper subgroup \(M<T\) has at most \(1+[T:M](|M|-1)<|T|\) elements.

We apply Scott's lemma \cite[p.~328]{Scott} to \(D\). A subgroup that projects onto every factor of a product of nonabelian simple groups is a product of diagonal subgroups on disjoint blocks of coordinates. Within each block, the coordinate values are determined by one of them through isomorphisms. Indeed, by Goursat's lemma \cite[Theorem~2.1]{BSZ}, the projection onto two coordinates is either the full product or the graph of an isomorphism. The graph projections correspond to equality of the kernels of the coordinate maps, so they define equivalence classes. To see that the projection onto one coordinate from each equivalence class is the full product, apply induction and Goursat's lemma. If the last projection were linked to the product of the preceding factors, there would be a surjection from that product onto the last simple factor. Such a surjection factors through an individual factor, because the images of the factors are commuting normal subgroups of the simple image. This would give a graph projection between two of the chosen representatives.

The partition for \(D\) is invariant under translations by \(V\), since \(H\) normalizes \(D\) and maps onto \(V\). Let \(U(H)\) be the block containing \(0\). Translation by any element of \(U(H)\) preserves this block, so \(U(H)\) is an additive subgroup of \(V\). Since \(\ell\) is prime, \(U(H)\) is a subspace. The other blocks are its cosets. We have \(U(H)\ne0\), since singleton blocks would give \(D=N\) and then \(H=G\).

Suppose that \(U(H)\) contains a forbidden line with chosen vector \(q\). Since \(H\) meets the class of \((f_q,q)\), choose \(h\in H\) in that class. Conjugation by \(h\) preserves every diagonal block of \(D\), since the translation \(q\) preserves the cosets of \(U(H)\). Let \(\beta\) be the automorphism induced on one diagonal block. To compute \((f_q,q)^\ell\), multiply the values of \(f_q\) along each translation cycle. Because exactly one value on each cycle is \(s\), this power is the constant base function with value \(s\). Thus every coordinate of \(h^\ell\) has order \(\ell\). On the chosen block, projection onto any coordinate identifies \(\beta^\ell\) with conjugation by that coordinate of \(h^\ell\). This inner automorphism has order \(\ell\) because \(T\) is centerless. Hence \(\beta\) has order \(\ell^2\), a contradiction.

For each nonzero subspace \(U\) containing no forbidden line, define
\[
H_U=\{(f,t)\in G:f\text{ is constant on each coset of }U\}.
\]
Then \(H_U\) is a proper subgroup and \(U(H_U)=U\) once we verify \(H_U\in P\). This subgroup meets the classes in (1) and (2) by taking constant functions and the elements \((1,e_j)\). For (3), the vector \(q\) is outside \(U\). Translation by \(q\) on \(V/U\) has cycles of length \(\ell\). Choose one coset in each cycle and set \(f=s\) on those cosets and \(f=1\) on the others. On each translation cycle of \(q\) in \(V\), both \(f\) and \(f_q\) take the value \(s\) exactly once. They give conjugate elements \((f,q)\) and \((f_q,q)\). To see this, conjugation by a base function \(b\) replaces \(f_q(v)\) by \(b(v)f_q(v)b(v-q)^{-1}\). Taking \(b\) with values in \(\langle s\rangle\), the equations on each cycle are consistent because the sums of the exponents of \(s\) agree.

Let \(F(K)\) be the poset of nonempty faces of \(K\). We now have maps
\[
 P\longrightarrow F(K),\quad H\longmapsto\supp(U(H)),
 \qquad
 F(K)\longrightarrow P,\quad F\longmapsto H_{\langle e_j:j\in F\rangle}.
\]
Both maps reverse inclusion. For the first map, if \(H_1\leq H_2\), then equality of two coordinate kernels on \(H_2\cap N\) implies equality of their restrictions to \(H_1\cap N\). Thus \(U(H_1)\supseteq U(H_2)\). The second map reverses inclusion because enlarging \(U\) imposes more equalities on the coordinates. Their composite on \(F(K)\) is the identity. An inclusion-reversing map still sends a chain to a chain, so each map induces a simplicial map on order complexes. The composite is the identity on the barycentric subdivision of \(K\). Thus that subdivision is a retract of \(\Delta(P)\). Joining these maps with \(S^{r-2}\) and applying Lemma~\ref{lem:reduction} gives the homotopy retract \(\Sigma^{r-1}|K|\) of \(|\Delta_R|\).
\end{proof}

\begin{example}\label{ex:rp2}
Let \(K\) be the simplicial complex with the maximal faces listed below. We write \(ijk\) for \(\{i,j,k\}\).
\[
123,\ 124,\ 135,\ 146,\ 156,\ 236,\ 245,\ 256,\ 345,\ 346.
\]
The complex contains every edge, and each vertex link is a five-cycle. It is a closed connected surface with Euler characteristic \(6-15+10=1\), hence is a triangulation of \(\RP^2\). Every minimal nonface is a triangle. Indeed, a four-element set with all of its triangles present would give a triangle in a vertex link. We can therefore take \(\ell=3\) and \(T=A_5\) in Theorem~\ref{thm:realization}.

In \(\F_3^6\), the number of vectors with support a face of \(K\), including the zero vector, is
\[
1+6\cdot 2+15\cdot 2^2+10\cdot 2^3=153.
\]
There are \((3^6-153)/2=288\) forbidden lines. The \(4\) nonidentity conjugacy classes of \(A_5\), the \(6\) basis vectors and the \(288\) forbidden lines give \(r=298\). These classes are distinct. The last two types have different images in \(V\) because the basis vectors have one-element supports, and the forbidden lines do not. Since \(H_1(K;\Z)\cong \Z/2\Z\), the theorem gives a direct summand \(\Z/2\Z\) in \(\widetilde H_{298}(\Delta_R;\Z)\).
\end{example}

\begin{corollary}\label{cor:abelian}
For each finite abelian group \(A\), there exist a finite conjugation quandle \(R\) and an integer \(j\) such that \(A\) is isomorphic to a direct summand of \(\widetilde H_j(\Delta_R;\Z)\).
\end{corollary}

\begin{proof}
If \(A=0\), the assertion is immediate. Otherwise, write \(A\) as a direct sum of finite cyclic groups. For each cyclic summand of order \(d\), take the mapping cone of a piecewise linear map \(S^1\to S^1\) of degree \(d\). Triangulate the mapping cylinder and attach a cone along its copy of the domain. This gives a finite simplicial model of the mapping cone. Its first integral homology is \(\Z/d\Z\) by the degree-\(d\) attaching map \cite[Chapter~2]{Hatcher}. Wedge the resulting finite complexes at vertices and apply Theorem~\ref{thm:realization}.
\end{proof}

\section{Reductions for finite acyclicity}\label{sec:acyclic}

We return to the acyclicity question in \cite[Section~5]{HSW}. Does there exist a finite rack \(R\) such that \(\Delta_R\) is nonempty and \(\widetilde H_i(\Delta_R;\Z)=0\) for each \(i\geq0\)? For a finite group \(G\), let \(\Rad(G)\) denote its largest solvable normal subgroup. We show that one may take \(\Rad(G)=1\) in the subgroup description of Lemma~\ref{lem:reduction}.

\begin{lemma}\label{lem:quotient}
Let \(C_1,\dots,C_r\) be conjugacy classes generating a finite group \(G\), and assume that \(N\) is an abelian minimal normal subgroup of \(G\). Write \(\pi:G\to G/N\) for the quotient map and let \(\overline C_i=\pi(C_i)\). Put \(P=\PP(G;C_1,\dots,C_r)\) and \(\overline P=\PP(G/N;\overline C_1,\dots,\overline C_r)\). Then \(\Delta(P)\) is acyclic if and only if \(\Delta(\overline P)\) is acyclic.
\end{lemma}

\begin{proof}
Let \(P_0\) be the subposet of \(P\) consisting of those \(H\) with \(HN<G\). Projection and inverse image give maps \(P_0\to\overline P\) and \(\overline P\to P_0\). The composite on \(\overline P\) is the identity. The composite on \(P_0\) sends each \(H\) to \(HN\), which contains \(H\). Hence \(\Delta(P_0)\simeq\Delta(\overline P)\).

If \(H\in P\setminus P_0\), then \(HN=G\). Since \(N\) is abelian, \(H\cap N\) is normalized by both \(H\) and \(N\). Minimality of \(N\) and properness of \(H\) give \(H\cap N=1\). Thus \(H\) is a complement to \(N\). The complements in \(P\) are incomparable. For each such \(H\), projection identifies its strict lower interval in \(P\) with \(\overline P\). Indeed, \(\pi|_H\) is an isomorphism, and \(H\cap C_i\) maps onto \(\overline C_i\) for each \(i\), because \(H\) meets \(C_i\) and maps onto \(G/N\). Therefore the inclusion of the lower interval into \(P_0\) is a homotopy equivalence on order complexes.

Let \(c\) be the number of complements in \(P\). We cannot have \(c=1\), since that complement would be normal in \(G\) by conjugation invariance, and meeting every \(C_i\) would then imply that it contains the generating union. If \(c=0\), then \(P=P_0\). If \(c\geq 2\) and \(\overline P\) is nonempty, then \(\Delta(P)\) is obtained by attaching \(c\) cones to \(\Delta(P_0)\) along subcomplexes whose inclusions are homotopy equivalences. The first cone together with \(\Delta(P_0)\) is contractible, since its attaching inclusion is a homotopy equivalence. Collapsing this contractible subcomplex to a point gives a homotopy equivalence \cite[Proposition~0.17]{Hatcher}. Each of the remaining cones then has its base identified to that point and gives a suspension of \(\Delta(\overline P)\). Thus \(\Delta(P)\) is homotopy equivalent to a wedge of \(c-1\) such suspensions. If \(\overline P\) is empty and \(c\geq2\), then \(\Delta(P)\) consists of \(c\) isolated points. These cases prove the assertion.
\end{proof}

\begin{theorem}\label{thm:acyclic}
The following statements are equivalent.
\begin{enumerate}
\item There is a finite rack \(R\) for which \(\Delta_R\) is acyclic.
\item There is a finite conjugation quandle \(R\) which is a generating union of conjugacy classes in a finite group \(G\) with \(\Rad(G)=1\), and for which \(\Delta_R\) is acyclic.
\item There is a finite conjugation quandle \(R\) as in \textup{(2)} for which \(\Delta_R\) is contractible.
\end{enumerate}
\end{theorem}

\begin{proof}
Assume (1). Lemma~\ref{lem:reduction} gives a nonempty acyclic class-meeting subgroup complex for \(G=\Inn(R)\). Indeed, if the subgroup poset is empty, the rack complex has the homotopy type of a sphere under the augmented convention. For a nonempty subgroup complex the join in Lemma~\ref{lem:reduction} shifts reduced homology.

Choose a chief series from \(1\) to \(\Rad(G)\). Each factor is abelian. Apply Lemma~\ref{lem:quotient} successively to the corresponding quotients. This gives an acyclic class-meeting subgroup complex in \(G/\Rad(G)\), which has trivial solvable radical. The projected classes still generate the quotient. Discard repetitions and let \(R'\) be the union of the remaining classes with the conjugation operation. Repetitions do not affect the subgroup poset. Lemma~\ref{lem:reduction} now shows that \(\Delta_{R'}\) is acyclic. Thus (1) implies (2).

Assume (2), and put \(R_0=R\setminus\{1\}\). The class-meeting subgroup poset is unchanged on removing or adjoining the identity class. The set \(R_0\) still generates \(G\), and Lemma~\ref{lem:reduction} shows that \(\Delta_{R_0}\) is acyclic and that \(\Delta_{R_0\cup\{1\}}\simeq\Sigma\Delta_{R_0}\). The suspension is simply connected and acyclic, hence contractible by the homology form of Whitehead's theorem \cite[Corollary~4.33]{Hatcher}. Thus (2) implies (3), and (3) implies (1).
\end{proof}

\section*{Use of generative AI}
Under Amir Saki's direction and supervision, OpenAI's ChatGPT was used in parts of the development and verification of the arguments and in the preparation of this manuscript. The authors take responsibility for the contents of the manuscript.

\end{document}